\documentclass[11pt]{article}

\usepackage{amsfonts,amssymb}
\usepackage{a4wide}
\usepackage[utf8]{inputenc}

\usepackage{etoolbox}

\makeatletter
\patchcmd{\maketitle}{\@fnsymbol}{\@alph}{}{}  
\makeatother

\newtheorem{theorem}{Theorem}

\newtheorem{lemma}{Lemma}
\newtheorem{proposition}{Proposition}
\newtheorem{innercustomthm}{Theorem}
\newenvironment{mythm}[1]
  {\renewcommand\theinnercustomthm{#1}\innercustomthm}
  {\endinnercustomthm}

\usepackage{hyperref}

\title{On Thurston's geometrical space form problem: \\ on quasi space forms}

\author{S. Haesen\thanks{Department of Teacher Education, Thomas More UC, campus Vorselaar, Belgium and Department of Mathematics, University of Hasselt, Belgium, e-mail address: stefan.haesen@thomasmore.be}%
  \and and
  \and
  M. Petrovi\'{c}-Torga\v{s}ev \thanks{State University of Novi Pazar, Vuka Karadzica 9, 36300 Novi Pazar, Serbia.}%
  \and and
  \and
  L. Verstraelen \thanks{Section of Geometry, University of Leuven, Belgium and CiT, De Haan, Belgium, e-mail address: leopold.verstraelen@kuleuven.be}%
  }

\date{Dedicated to Distinguished Professor Bang-Yen Chen \\ at the occasion of his 80th anniversary}

\begin{document}

\maketitle

\begin{abstract}
A proposal is made for what may well be the most elementary Riemannian spaces which are homogeneous but not isotropic. In other words: a proposal is made for what may well be the \textit{the nicest symmetric spaces beyond the real space forms}, that is, \textit{beyond the Riemannian spaces which are homogeneous and isotropic}. The above qualification of `'nicest symmetric spaces'' finds a justification in that, \textit{together with the real space forms, these spaces are most natural with respect to the importance in human vision of our ability to readily recognise conformal things} and in that \textit{ these spaces are most natural with respect to what in Weyl's view is symmetry in Riemannian geometry}.

Following his suggestion to remove the real space forms' isotropy condition, \textit{the quasi space forms thus introduced do offer a metrical, local geometrical solution to the geometrical space form problem as posed by Thurston in his 1979 Princeton Lecture Notes} on `'The Geometry and Topology of 3-manifolds''. Roughly speaking, quasi space forms are the Riemannian manifolds of dimension greater than or equal to 3, which are not real space forms but which admit two orthogonally complementary distributions such that at all points all the 2-planes that in the tangent spaces there are situated in a same position relative to these distributions do have the same sectional curvatures. 
\end{abstract}


\section*{Introduction}

This article is to be classified in the category `'\textit{Local Riemannian Geometry}'' within the field of `'\textit{Classical Differential Geometry}''. It consists of the following sections: (1) Definition of quasi space forms; (2) On the symmetric spaces of Deszcz; (2) Definition of quasi Einstein spaces; (4) The conformally Euclidean Deszcz symmetric spaces.

As usual, this paper ends with a \textit{list of references} of which, for what comes next, in particular \textit{the lecture of the two following general articles might be not without interest}: S. Haesen and L. Verstraelen, `'\textit{Natural Intrinsic Geometrical Symmetries}'', (in SIGMA 5 - Special Issue `'Elie Cartan and Differential Geometry''; 2009), and, L. Verstraelen, `'\textit{Submanifold theory - A contemplation of submanifolds}'', (in the AMS series Contemporary Mathematics, book 756, from the `'AMS Special Session on \textit{Geometry of Submanifolds} in honor of Bang-Yen Chen''; 2020). And, as \textit{basic reference for the fundamentals of differential geometry in general and for some classical theorems concerning the contents of the present paper in particular}, we recommend to consult \textit{K\"{u}hnel's Lecture Notes}: W. K\"{u}hnel, `\textit{'Differentialgeometrie. Kurven - Fl{\"a}chen - Mannigfaltigkeiten}'', (4. Auflage, Vieweg Studium; 2008), in English translation: `'\textit{Differential Geometry: Curves - Surfaces - Manifolds}'', (AMS Student Mathematical Library, volume 16; 2005).


\section{Definition of quasi space forms}

A Riemannian manifold $(M^{n},g)$ of dimension $n\geq 3$ is said to be \textit{a real space form} if it is \textit{a space of constant curvature}, for short, \textit{a CC space}, that is, if \textit{for all its tangent 2-planes $\tilde{\pi}$ at all of its points its Riemannian sectional curvatures $K(\tilde{\pi})$ are the same}, that is, if \textit{for all $\tilde{\pi}$ at all points $K(\tilde{\pi})=c$ for one and the same real number $c=0$, $c>0$ or $c<0$}, and \textit{such spaces are denoted by} $M^{n}(c)$. For basic information on real space forms, see \cite{kuhnel}.

A Riemannian manifold $(M^{n},g)$ of dimension $n\geq 3$ will be said to be \textit{a quasi space form} or \textit{a space of quasi constant curvature}, for short, \textit{a QCC space}, or more specifically, \textit{for some q}, ($1\leq q\leq k$; according to $n$ being even, $n=2k$, or $n$ being odd, $n=2k+1$), $(M^{n},g)$ is said to be \textit{a (q) quasi space form} or \textit{ a space of (q) quasi constant curvature}, for short, \textit{a (q) QCC space, if $(M^{n},g)$ admits two orthogonally complementary distributions $D$ and $D^{\perp}$, of  fixed dimensions $q$ and $q^{\perp}=n-q\geq q$ respectively, and if its Riemannian sectional curvatures $K(\tilde{\pi})$ do depend on the positions of the tangent 2-planes $\tilde{\pi}$ with respect to the direct sum decomposition $TM = D\oplus D^{\perp}$ in the way that hereafter will be stated separately for the cases $q=1$ and $q>1$}. 

{\noindent}\textit{In case $q=1$}: for all orthonormal $X\in D$ and $X^{\perp}, Y^{\perp}\in D^{\perp}$ and for all $\theta\in \left[0,\pi/2\right]$, \textit{the curvature $K(\tilde{\pi})$ for all planes $\tilde{\pi}=(X\, \cos\theta+X^{\perp}\sin\theta)\wedge Y^{\perp}$ is given by}

\begin{equation}
K(\tilde{\pi}) = \overline{K}\cos^{2}\theta + K^{\perp}\sin^{2}\theta\ , \label{eq1}
\end{equation}

{\noindent}\textit{whereby $\overline{K}$ and $K^{\perp}$ are everywhere distinct real functions on $M^{n}$}.

{\noindent}\textit{In case $q>1$}: for all orthonormal $X,Y\in D$ and $X^{\perp},Y^{\perp}\in D^{\perp}$ and for all $\theta,\varphi\in \left[0,\pi/2\right]$, \textit{the curvature $K(\tilde{\pi})$ for all planes $\tilde{\pi} = (X\, \cos\theta+X^{\perp}\sin\theta)\wedge (Y\, \cos\varphi+Y^{\perp}\sin\varphi)$ is given by}

\begin{equation}
K(\tilde{\pi}) = K\cos^{2}\theta\cos^{2}\varphi + K^{\perp}\sin^{2}\theta\sin^{2}\varphi + \overline{K}(\cos^{2}\theta\sin^{2}\varphi+\sin^{2}\theta\cos^{2}\varphi)\ , \label{eq2}
\end{equation} 

{\noindent}\textit{whereby $K$ and $K^{\perp}$ are everywhere distinct real functions on $M^{n}$ and $\overline{K}=(K+K^{\perp})/2$}.

\textit{Concerning (1) QCC spaces next follow some specifications of the so to say three different  possible kinds of tangent 2-planes $\tilde{\pi}$, according to their positions with respect to the distributions $D$ and $D^{\perp}$}. In the present situation, $D$ is 1-dimensional, say, $D$ is generated by a  unit vector $X$, $D=\left[X\right]$, and $D^{\perp}$ is of dimension $q^{\perp}=n-1\geq 2$, and $X^{\perp},Y^{\perp}$ denote arbitrary orthonormal vectors in $D^{\perp}$. For $\theta=\pi/2$, the planes $\tilde{\pi}=(X\, \cos\theta+X^{\perp}\sin\theta)\wedge Y^{\perp} = X^{\perp}\wedge Y^{\perp}$ will be denoted by $\pi^{\perp}$: the planes $\pi^{\perp}$ lie in $D^{\perp}$ and are perpendicular to $D$, or, still, $\angle(\pi^{\perp},D)=\pi/2$ and $\angle(\pi^{\perp},D^{\perp})=0$; and (\ref{eq1}) states that all planes $\pi^{\perp}$ do have the same sectional curvature $K(\pi^{\perp})=K^{\perp}$. For $\theta=0$, the planes $\tilde{\pi}=(X\cos\theta+X^{\perp}\sin\theta)\wedge Y^{\perp} = X\wedge Y^{\perp}$ will be denoted by $\overline{\pi}$: the planes $\overline{\pi}$ contain the line $\left[ X\right] = D$ and cut $D^{\perp}$, as well as project onto $D^{\perp}$, in the line $\left[Y^{\perp}\right]$, or still, $\angle(\overline{\pi},D) = 0$ and $\angle(\overline{\pi},D^{\perp}) = \pi/2$; and (\ref{eq1}) states that all planes $\overline{\pi}$ through $D$ do have the same sectional curvature $K(\overline{\pi})=\overline{K}\neq K^{\perp}$. For $\theta\neq 0$ and $\theta\neq \pi/2$, the vectors $\tilde{X}=X\cos\theta+X^{\perp}\sin\theta$ do properly position in between $D$ and $D^{\perp}$, and the corresponding planes $\tilde{\pi}=(X\cos\theta+X^{\perp}\sin\theta)\wedge Y^{\perp}$ will be denoted by $\pi_{\theta}$: planes $\pi_{\theta}$ project onto $D$ in the line $\left[X\right]=D$ with which line they make an angle $\theta$ and planes $\pi_{\theta}$ project onto $D^{\perp}$ in the plane $X^{\perp}\wedge Y^{\perp}$ with which plane they make an angle $\theta^{\perp}=(\pi/2)-\theta$; (we recall that \textit{the angle $\psi$ between two 2-planes $\pi^{\alpha}=\tilde{A}_{1}\wedge\tilde{A}_{2}$ and $\pi^{\beta} = \tilde{B}_{1}\wedge\tilde{B}_{2}$ which are given by orthonormal vectors $\tilde{A}_{1},\tilde{A}_{2}$ and $\tilde{B}_{1},\tilde{B}_{2}$ respectively, in an arbitrary dimensional Euclidean space, $\psi=\angle(\pi^{\alpha},\pi^{\beta})$, is determined by $\cos^{2}\psi=(\mbox{det}\, M)^{2}$, whereby $M$ is the $2\times 2$ matrix with elements $M_{11}=g(\tilde{A}_{1},\tilde{B}_{1})$, $M_{12}=g(\tilde{A}_{1},\tilde{B}_{2})$, $M_{21}=g(\tilde{A}_{2},\tilde{B}_{1})$ and $M_{22}=g(\tilde{A}_{2},\tilde{B}_{2})$}, \cite{risteki}); and (\ref{eq1}) states that all planes $\pi_{\theta}$, that is, all tangent 2-planes that make an angle $\theta$ with $D$, do have the same sectional curvature $K(\pi_{\theta})=\overline{K}\cos^{2}\theta+K^{\perp}\sin^{2}\theta$.

\textit{Concerning (q)QCC spaces for $q>1$ next follow some specifications of the so to say six different possible kinds of tangent 2-planes $\tilde{\pi}$, according to their positions with respect to the distributions $D$ and $D^{\perp}$.} In the present situation, $D$ and $D^{\perp}$ have dimensions $q$ and $q^{\perp}=n-q$ respectively, which both are $\geq 2$. Every tangent 2-plane $\tilde{\pi}$ has an orthonormal basis formed by vectors $\tilde{X} = X\cos\theta+X^{\perp}\sin\theta$ and $\tilde{Y}=Y\cos\varphi+Y^{\perp}\sin\varphi$, for some angles $\theta,\varphi\in \left[0,\pi/2\right]$ and for some orthonormal $X,Y\in D$ and for some orthonormal $X^{\perp},Y^{\perp}\in D^{\perp}$. For $\theta=\varphi=\pi/2$, the planes $\tilde{\pi}=(X\cos\theta+X^{\perp}\sin\theta)\wedge (Y\cos\varphi+Y^{\perp}\sin\varphi) = X^{\perp}\wedge Y^{\perp}$ will be denoted by $\pi^{\perp}$: the planes $\pi^{\perp}$ lie in $D^{\perp}$ and are perpendicular to $D$, or, still, $\angle(\pi^{\perp},D)=\pi/2$ and $\angle(\pi^{\perp},D^{\perp})=0$; and (\ref{eq2}) states that all planes $\pi^{\perp}$ do have the same sectional curvature $K(\pi^{\perp})=K^{\perp}$. For $\theta=\varphi=0$, the planes $\tilde{\pi}=(X\cos\theta+X^{\perp}\sin\theta)\wedge (Y\cos\varphi+Y^{\perp}\sin\varphi)=X\wedge Y$ will be denoted by $\pi$: the planes $\pi$ lie in $D$ and are perpendicular to $D^{\perp}$, or, still, $\angle(\pi,D)=0$ and $\angle(\pi,D^{\perp})=\pi/2$; and (\ref{eq2}) states that all planes $\pi$ do have the same sectional curvature $K(\pi)=K\neq K^{\perp}$. For (i) $\theta=0$ and $\varphi=\pi/2$, and, similarly, for (ii) $\theta=\pi/2$ and $\varphi=0$, the planes $\tilde{\pi}=(X\cos\theta+X^{\perp}\sin\theta)\wedge (Y\cos\varphi+Y^{\perp}\sin\varphi)$ are either given by (i) $\tilde{\pi}=X\wedge Y^{\perp}$ or by (ii) $\tilde{\pi}=X^{\perp}\wedge Y$, and such planes will be denoted by $\overline{\pi}$: both onto $D$ and onto $D^{\perp}$, planes $\overline{\pi}$ project in a line, and $\angle(\overline{\pi},D)=\angle(\overline{\pi},D^{\perp})=\pi/2$; and (\ref{eq2}) states that all planes $\overline{\pi}$ do have the same sectional curvature $\overline{K}=(K+K^{\perp})/2$. Further, we consider the situations (i) $0<\theta<\pi/2$ and $\varphi=\pi/2$, and, (ii) $\theta=\pi/2$ and $0<\varphi<\pi/2$; geometrically they amount to the same and so we will confine here only to the case (i), whereby the vectors $\tilde{X}=X\cos\theta+X^{\perp}\sin\theta$ do properly position in between $D$ and $D^{\perp}$, and the planes $\tilde{\pi}=(X\cos\theta+X^{\perp}\sin\theta)\wedge (Y\cos\varphi+Y^{\perp}\sin\varphi) = (X\cos\theta+X^{\perp}\sin\theta)\wedge Y^{\perp}$ will be denoted by $\pi_{\theta}$: planes $\pi_{\theta}$ project onto $D$ in the line $\left[X\right]$ with which they make an angle $\theta$ and planes $\pi_{\theta}$ project onto $D^{\perp}$ in the plane $X^{\perp}\wedge Y^{\perp}$ with which they make an angle $\theta^{\perp} = \pi/2-\theta$; and (\ref{eq2}) states that all planes $\pi_{\theta}$ do have the same sectional curvature $K_{\theta}=\overline{K}\cos^{2}\theta+K^{\perp}\sin^{2}\theta$. Still further, we consider the situations (i') $0<\theta<\pi/2$ and $\varphi=0$, and, (ii') $\theta=0$ and $0<\varphi<\pi/2$; geometrically they amount to the same and so we will confine here only to do the case (i'), whereby the vectors $\tilde{X}=X\cos\theta+ X^{\perp}\sin\theta$ do properly position in between $D$ and $D^{\perp}$, and the planes $\tilde{\pi}=(X\cos\theta+X^{\perp}\sin\theta)\wedge (Y\cos\varphi+Y^{\perp}\sin\varphi)=(X\cos\theta+X^{\perp}\sin\theta)\wedge Y$ will be denoted by $\pi^{\perp}_{\theta^{\perp}}$: planes $\pi^{\perp}_{\theta^{\perp}}$ project onto $D^{\perp}$ in the line $\left[X^{\perp}\right]$ with which they make an angle $\theta^{\perp}$ and planes $\pi^{\perp}_{\theta^{\perp}}$ project onto $D$ in the plane $X\wedge Y$ with which they make an angle $\theta$; and (\ref{eq2}) states that all planes $\pi^{\perp}_{\theta^{\perp}}$ do have the same sectional curvature $K(\pi^{\perp}_{\theta^{\perp}}) = \overline{K}\sin^{2}\theta+K\cos^{2}\theta (=\overline{K}\cos^{2}\theta^{\perp}+K\sin^{2}\theta^{\perp})$. Finally, for $0<\theta,\varphi<\pi/2$, the planes $\tilde{\pi}=(X\cos\theta+X^{\perp}\sin\theta)\wedge (Y\cos\varphi+Y^{\perp}\sin\varphi)$, whereby both vectors $\tilde{X}=X\cos\theta+X^{\perp}\sin\theta$ and $\tilde{Y}=Y\cos\varphi+Y^{\perp}\sin\varphi$ do properly position in between $D$ and $D^{\perp}$ will be denoted by $\pi_{\theta,\varphi}$ or by $\pi^{\perp}_{\theta^{\perp},\varphi^{\perp}}$, whereby $\theta^{\perp}=\pi/2-\theta$ and $\varphi^{\perp}=\pi/2-\varphi$, (using one or other notation as it may feel to match better with the occasion at hand or just like that): planes $\pi_{\theta,\varphi}=\pi^{\perp}_{\theta^{\perp},\varphi^{\perp}}$ project onto $D$ in the plane $X\wedge Y$ with which they make an angle $\arccos(\cos\theta\cos\varphi)$ and planes $\pi^{\perp}_{\theta^{\perp},\varphi^{\perp}}=\pi_{\theta,\varphi}$ project onto $D^{\perp}$ in the plane $X^{\perp}\wedge Y^{\perp}$ with which they make an angle $\arccos(\cos\theta^{\perp}\cos\varphi^{\perp})$; and (\ref{eq2}) states that all planes $\pi_{\theta,\varphi}=\pi^{\perp}_{\theta^{\perp},\varphi^{\perp}}$ do have the same sectional curvature 

\[ K(\pi_{\theta,\varphi}) = K(\pi^{\perp}_{\theta^{\perp},\varphi^{\perp}}) = K\cos^{2}\theta\cos^{2}\varphi + K^{\perp}\cos^{2}\theta^{\perp}\cos^{2}\varphi^{\perp}+\overline{K} (\cos^{2}\theta\cos^{2}\varphi^{\perp} + \cos^{2}\theta^{\perp}\cos^{2}\varphi)\ . \]

Corresponding to the notation $M^{n}(c)$ for the $n$-dimensional CC spaces of curvature $c$, in view of the above definition, \textit{the $n$-dimensional QCC spaces of curvatures $\overline{K}$ and $K^{\perp}$ and $K=2\overline{K}-K^{\perp}\neq K^{\perp}$ might well be denoted by $M^{n}(\overline{K},K^{\perp})$}. Also, \textit{(q) quasi constant curvature spaces or (q) QCC spaces may sometimes simply be called quasi constant curvature spaces or QCC spaces}, and, then the actual dimensions $q\geq 1$ and $q^{\perp}=n-q\geq q$ of their distributions $D$ and $D^{\perp}$ will only be specified when needed.

As far as we know, first purposefull studies of the (1)QCC spaces were done in the early 1970ties by B.-Y. Chen and K. Yano and C.-S. Houh \cite{chen1,chen2,yano1}, and, from the later studies that we know of, here we would like to mention in particular those by V. Boju and M. Popescu \cite{boju} and by G. Ganchev and V. Mihova \cite{ganchev1,ganchev2}.


\section{On the symmetric spaces of Deszcz}

\textit{A geometrical symmetry of a Riemannian space} $(M^{n},g)$ concerns the invariance of some measure of `'some geometrical beings that live on the manifold $M^{n}$'', under the performance of some kind of transformations \cite{trans}. The transformations of Riemannian spaces $(M^{n},g)$ that we will consider hereafter are \textit{the parallel transports fully around infinitesimal co-ordinate parallelograms on the manifold $M^{n}$}; (as described in some detail in \cite{trans2}, these transformations do show utmost respect both for the differential structure and for the metrical structure of Riemannian spaces). And then various corresponding geometrical symmetries may be studied on spaces $(M^{n},g)$ depending on which measures are taken of which beings.

In \cite{schouten} Schouten showed that \textit{at all points $p$ of all Riemannian manifolds, the changes of the directions of their tangent vectors under their parallel transport fully around the infinitesimal co-ordinate parallelograms cornered at $p$ are measured by the Riemann-Christoffel curvature tensor $R$ of these spaces}, and thus he obtained the following.

\vspace{2mm}

\begin{mythm}{A}
The Riemannian spaces for which all tangent vectors at all points remain invariant under the parallel transport fully around all infinitesimal co-ordinate parallelograms are the locally Euclidean spaces, or, still, the locally flat spaces.
\end{mythm}

\vspace{2mm}

{\noindent}\textit{The locally Euclidean spaces} are the Riemannian spaces $(M^{n},g)$ with \textit{vanishing curvature tensor R}, $R=0$, or, still, it are the Riemannian spaces $(M^{n},g)$ with \textit{vanishing sectional curvatures $K(p,\tilde{\pi})$}, $K(p,\tilde{\pi})=0$, for all points $p\in M^{n}$ and for all 2-planes $\tilde{\pi}\subset T_{p}M^{n}$; it are \textit{the CC spaces $M^{n}(0)$}. And, next follows a well known \textit{Theorem of Beltrami}, as reference for which we suggest Vladimir Matveev's article \cite{matveev}.

\vspace{2mm}

\begin{mythm}{B}
The real space forms constitute the projective class of the locally Euclidean spaces: by applying geodesic transformations to CC spaces $M^{n}(0)$ one obtains CC spaces $M^{n}(c)$ of constant curvature $c=0$ or $c>0$ or $c<0$ and the class of all CC spaces $M^{n}(c)$ is closed under geodesic transformations.
\end{mythm}

\vspace{2mm}

\textit{The Riemannian sectional curvatures} $K(p,\tilde{\pi})$ at points $p\in M^{n}$ and for tangent 2-planes $\tilde{\pi}\subset T_{p}M^{n}$ \textit{are the main metrical invariants of Riemannian manifolds} $(M^{n},g)$, (cfr. \cite{berger1,berger2}). By the parallel transports of any 2-plane $\tilde{\pi}\subset T_{p}M^{n}$ at any point $p\in M^{n}$ fully around the infinitesimal co-ordinate parallelograms cornered at $p$ one obtains 2-planes $\tilde{\pi}^{\star}\subset T_{p}M^{n}$ whereby, in general, $\tilde{\pi}^{\star}\neq \tilde{\pi}$. In \cite{haesen2} it was shown that \textit{at all points $p$ of all Riemannian manifolds and for all tangent 2-planes $\tilde{\pi}$ there, the changes of their sectional curvatures under such parallel transports}, that is, the changes $K(p,\tilde{\pi}^{\star}) - K(p,\tilde{\pi})$, \textit{are measured by the curvature tensor $R\cdot R$ of these spaces}, (in $R\cdot R$, \textit{the second $R$} stands for the Riemann-Christoffel \textit{curvature tensor} and \textit{the first $R$} stands for \textit{the curvature operator which}, by the meaning of the dot $\cdot$ \textit{acts as a derivation on the curvature tensor} $R$), which yields the following. 

\vspace{2mm}

\begin{mythm}{C}
The Riemannian spaces for which all Riemannian sectional curvatures remain invariant under the parallel transports fully around all infinitesimal co-ordinate parallelograms are the Szab\'{o} symmetric spaces, or, still, the semi symmetric spaces.
\end{mythm}

\vspace{2mm}

{\noindent}\textit{The Szab\'{o} symmetric spaces} are the Riemannian spaces $(M^{n},g)$ with \textit{vanishing curvature tensor $R\cdot R$}, $R\cdot R=0$; such Riemannian spaces are also said to be semi symmetric and these spaces were classified by Zoltan Szab\'{o} \cite{szabo1,szabo2}. \textit{The locally symmetric spaces}, or, still, \textit{the Cartan symmetric spaces}, are the Riemannian spaces $(M^{n},g)$ for which \textit{the curvature tensor $R$ is parallel}, $\nabla R=0$, whereby $\nabla$ denotes the \textit{Levi-Civita connection} of $(M^{n},g)$, or, still, by a \textit{Theorem of Levy} \cite{levy}, it are \textit{the Riemannian spaces for which all Riemannian sectional curvatures $K(p,\tilde{\pi})$ remain invariant under the parallel transport of $\tilde{\pi}$ along all infinitesimal geodesics emanating from $p$ on $M^{n}$}. These spaces were classified by \'{E}lie Cartan, who moreover could characterise them as \textit{the Riemannian spaces for which the local geodesic reflections in all their points are local isometries}. As far as we know, the curvature condition $R\cdot R=0$ did first appear in the studies of locally symmetric spaces by \'{E}. Cartan and by P.A. Shirokov, namely as the integrability condition of $\nabla R=0$, and for more information on Cartan symmetric spaces and on Szab\'{o} symmetric spaces, in particular we refer to \cite{mikes1,mikes2,boeckx,lumiste}.

\textit{Given on a Riemannian manifold $(M^{n},g)$ a point $p$ and given a tangent 2-plane $\tilde{\pi}$ at $p$, then the Riemannian sectional curvature $K(p,\tilde{\pi})$ is a corresponding isometrically invariant scalar quantity}, namely the Gauss curvature at $p$ of the surface formed in $M^{n}$ around $p$ by the geodesics of $(M^{n},g)$ through $p$ of which the tangent line at $p$ lies in $\tilde{\pi}$. And, as shown by \'{E}lie Cartan, \textit{the knowledge of all the  sectional curvatures $K(p,\tilde{\pi})$ of $(M^{n},g)$ is equivalent to the knowledge of the curvature tensor $R$ of $(M^{n},g)$}. Similarly, \textit{given on a Riemannian manifold $(M^{n},g)$ a point $p$ and given two tangent 2-planes $\tilde{\pi}_{1}$ and $\tilde{\pi}_{2}$ at $p$, then, as some corresponding isometrically invariant scalar quantity, a double sectional curvature $L(p,\tilde{\pi}_{1},\tilde{\pi}_{2})$ has been defined in \cite{haesen2}}, and \textit{the knowledge of all double sectional curvatures $L(p,\tilde{\pi}_{1},\tilde{\pi}_{2})$ of $(M^{n},g)$ is equivalent to the knowledge of the curvature tensor $R\cdot R$ of $(M^{n},g)$}.

{\noindent}\textit{The real space forms $M^{n}(c)$} essentially are defined as the Riemannian manifolds $(M^{n},g)$ for which at all points $p$ the sectional curvatures $K(p,\tilde{\pi})$ are independent of the 2-planes $\tilde{\pi}\subset T_{p}M^{n}$, or, still, for which \textit{at all points the Riemannian sectional curvature function is isotropic}, that is, equals the same value $c$ in all 2-dimensional tangent directions $\tilde{\pi}$. And, similarly, \textit{the Deszcz symmetric spaces $M^{n}((L))$} basically are defined as the Riemannian manifolds $(M^{n},g)$ for which at all points $p$ the double sectional curvatures $L(p,\tilde{\pi}_{1},\tilde{\pi}_{2})$ are independent of the two 2-planes $\tilde{\pi}_{1},\tilde{\pi}_{2}\subset T_{p}M^{n}$, or, still, for which \textit{at all points the double sectional curvature function is isotropic}, that is, $L(p,\tilde{\pi}_{1},\tilde{\pi}_{2})$ equals the same value $L(p)$ of some function $L:\, M^{n}\rightarrow \mathbb{R}$ for all tangent 2-planes $\tilde{\pi}_{1}$ and $\tilde{\pi}_{2}$ at $p$ for which $L(p,\tilde{\pi}_{1},\tilde{\pi}_{2})$ is well defined. When the sectional curvatures $K(p,\tilde{\pi})$ do not depend on the planes $\tilde{\pi}$ at $p$, then, by \textit{the lemma of Schur}, the sectional curvatures $K(p,\tilde{\pi})$ moreover do not depend on the points $p$, such that then $\forall\, p\in M^{n}\ \forall\, \tilde{\pi}\subset T_{p}M^{n}:\ K(p,\tilde{\pi})=c$ for some fixed real number $c$, (cfr. \cite{kuhnel}). However, in the situation that the double sectional curvatures $L(p,\tilde{\pi}_{1},\tilde{\pi}_{2})$ do not depend on the planes $\tilde{\pi}_{1}$ and $\tilde{\pi}_{2}$ at $p$, then there is no such kind of lemma: \textit{in general, the double sectional curvature of a Deszcz symmetric space is a non-constant function $L:\, M^{n}\rightarrow\mathbb{R}$}. And the very special particular \textit{Deszcz symmetric spaces for which $L$ is a constant function}, following Kowalski and Sekizawa \cite{kowalski1,kowalski2} \textit{are said to be Deszcz symmetric spaces of constant type}. The Szab\'{o} symmetric spaces are the Riemannian spaces with vanishing curvature tensor $R\cdot R$, $R\cdot R=0$, or, still, it are \textit{the Deszcz symmetric spaces $M^{n}((0))$ of constant type $0$}; $0$ is hereby the zero function on $M^{n}$, namely $0:\, M^{n}\rightarrow \mathbb{R}:\, p\mapsto 0(p)=0$. And, in analogy with the above Theorem of Beltrami, due to Sinyukov, Mike\v{s}, Venzi, Defever and Deszcz \cite{defever,mikes1,mikes2}, one has the following.

\vspace{2mm}

\begin{mythm}{D}
The Deszcz symmetric spaces $M^{n}((L))$ constitute the projective class of the Szab\'{o} symmetric spaces: by applying geodesic transformations to Szab\'{o} symmetric spaces $M^{n}((0))$ one obtains Deszcz symmetric spaces $M^{n}((L))$ of any double sectional curvature function $L:\, M^{n}\rightarrow\mathbb{R}$, and the class of all Deszcz symmetric spaces $M^{n}((L))$ is closed under geodesic mappings.
\end{mythm}

\vspace{2mm}

In the following, by $\tilde{X},\tilde{Y},\tilde{V},\tilde{W},\tilde{V}_{1},\tilde{V}_{2},\tilde{V}_{3},\tilde{V}_{4}$ will be denoted arbitrary tangent vector fields on a Riemannian manifold $(M^{n},g)$. \textit{The curvature operators} $R(\tilde{X},\tilde{Y})=\nabla_{\tilde{X}}\nabla_{\tilde{Y}}-\nabla_{\tilde{Y}}\nabla_{\tilde{X}}-\nabla_{\left[\tilde{X},\tilde{Y}\right]}$ do \textit{determine the $(0,4)$ Riemann-Christoffel curvature tensor $R$ as follows}: $R(\tilde{X},\tilde{Y},\tilde{V},\tilde{W}) = g(R(\tilde{X},\tilde{Y})\tilde{V},\tilde{W})$. Similarly, \textit{the metrical endomorphisms} $\tilde{X}\wedge_{g}\tilde{Y}$ defined by $(\tilde{X}\wedge_{g}\tilde{Y})\tilde{V} = g(\tilde{Y},\tilde{V})\tilde{X}-g(\tilde{X},\tilde{V})\tilde{Y}$ do \textit{determine the $(0,4)$ tensor} $G(\tilde{X},\tilde{Y},\tilde{V},\tilde{W}) = g((\tilde{X}\wedge_{g}\tilde{Y})\tilde{V},\tilde{W})$, \textit{which likely is the simplest $(0,4)$ generalised curvature tensor on $(M^{n},g)$}, (cfr. \cite{kuhnel}); (\textit{the metrical endomorphisms $\vec{x}\wedge_{g}\vec{y}$, whereby $\vec{x},\vec{y}\in T_{p}M^{n}$, do measure the changes of the tangent vectors $\vec{v}$ at the points $p$ under their infinitesimal rotations with respect to the tangent 2-planes $\vec{x}\wedge\vec{y}$ at $p$}; see \cite{haesen2, trans2}). And, \textit{the sectional curvatures $K(p,\tilde{\pi})$ for tangent 2-planes $\tilde{\pi} = \vec{v}\wedge\vec{w}$ at points $p$ are defined by} 

\[ K(p,\tilde{\pi}) = \frac{R(\vec{v},\vec{w},\vec{w},\vec{v})}{G(\vec{v},\vec{w},\vec{w},\vec{v})}\ , \]

{\noindent}and thus can be considered as kind of calibrations of the action of the curvature operators of Riemannian manifolds $(M^{n},g)$ on tangent vectors by means of the action of the metrical endomorphisms of these manifolds $(M^{n},g)$ on the same tangent vectors. On Riemannian manifolds $(M^{n},g)$, that is, for definite metrics $g$, the above definition holds perfectly well for all possible tangent 2-planes $\vec{v}\wedge\vec{w}$ at all points $p$, whereas, for instance \textit{on semi-Riemannian manifolds} $(M^{n},g^{\#})$, that is, \textit{on spaces with indefinite metrics} $g^{\#}$, \textit{the definition of the sectional curvature of course only works well for the non-degenerate 2-planes $\vec{v}\wedge\vec{w}$}, i.e., only works in case $G(\vec{v},\vec{w},\vec{w},\vec{v}) = g^{\#}(\vec{v},\vec{v})g^{\#}(\vec{w},\vec{w}) - g^{\#}(\vec{v},\vec{w})^{2}\neq 0$. A geometrical interpretation for the adapted definition of the null sectional curvature in a Lorentzian manifold was given in \cite{albujer}.

\textit{On a Riemannian manifold $(M^{n},g)$, by their action as a derivation on the $(0,4)$ curvature tensor $R$, the curvature operators do determine the $(0,6)$ curvature tensor $R\cdot R$ as follows}:

\begin{eqnarray*} 
(R\cdot R)(\tilde{V}_{1},\tilde{V}_{2},\tilde{V}_{3},\tilde{V}_{4};\tilde{X},\tilde{Y})  & = & -R(R(\tilde{X},\tilde{Y})\tilde{V}_{1},\tilde{V}_{2},\tilde{V}_{3},\tilde{V}_{4}) - R(\tilde{V}_{1},R(\tilde{X},\tilde{Y})\tilde{V}_{2},\tilde{V}_{3},\tilde{V}_{4}) \\ 
& &  -R(\tilde{V}_{1},\tilde{V}_{2},R(\tilde{X},\tilde{Y})\tilde{V}_{3},\tilde{V}_{4}) - R(\tilde{V}_{1},\tilde{V}_{2},\tilde{V}_{3},R(\tilde{X},\tilde{Y})\tilde{V}_{4})\ .
\end{eqnarray*}

{\noindent}Similarly, \textit{by their action as a derivation on the $(0,4)$ curvature tensor $R$, the metrical endomorphisms do determine the $(0,6)$ so-called Tachibana tensor $\wedge_{g}\cdot R$. The Tachibana tensor $\wedge_{g}\cdot R$ measures the changes of the sectional curvatures $K(p,\tilde{\pi}_{2})$ of all tangent 2-planes $\tilde{\pi}_{2}=\vec{v}\wedge\vec{w}$ at all points $p$ under the infinitesimal rotations of these planes $\tilde{\pi}_{2}$ with respect to tangent 2-planes $\tilde{\pi}_{1}=\vec{x}\wedge\vec{y}$ at $p$, \cite{haesen2}, and the Tachibana tensor of a Riemannian manifold $(M^{n},g)$, $n\geq 3$, vanishes, $\wedge_{g}\cdot R = 0$, if and only if $(M^{n},g)$ is a real space form $M^{n}(c)$}, \cite{eisenhart}. And, \textit{the double sectional curvatures $L(p,\tilde{\pi}_{1},\tilde{\pi}_{2})$ for tangent 2-planes $\tilde{\pi}_{1}=\vec{x}\wedge\vec{y}$ and $\tilde{\pi}_{2}=\vec{v}\wedge\vec{w}$ at points $p$ are defined by} 

\[ L(p,\tilde{\pi}_{1},\tilde{\pi}_{2}) = \frac{(R\cdot R)(\vec{v},\vec{w},\vec{w},\vec{v};\vec{x},\vec{y})}{(\wedge_{g}\cdot R)(\vec{v},\vec{w},\vec{w},\vec{v};\vec{x},\vec{y})}\ , \]

{\noindent}and thus can be considered as kind of calibrations of the action of the curvature operators on the Riemann-Christoffel curvature $R$ by the action of the metrical endomorphisms on this same tensor $R$. Such double sectional curvatures $L(p,\tilde{\pi}_{1},\tilde{\pi}_{2})$ however, of course, are only well defined in case $(\wedge_{g}\cdot R)(\vec{v},\vec{w},\vec{w},\vec{v};\vec{x},\vec{y})\neq 0$, \textit{in which case the tangent 2-planes $\tilde{\pi}_{1}$ and $\tilde{\pi}_{2}$ are said to be curvature-dependent}; this phenomenon corresponds to sectional curvatures $K(p,\tilde{\pi})$ on indefinite spaces $(M^{n},g^{\#})$ only being well defined for non-degenerate tangent 2-planes $\tilde{\pi}$. And, in particular, like \textit{a Riemannian manifold $(M^{n},g)$ is algebraically characterised to be a real space form $M^{n}(c)$ by the constant scalar valued proportionality of its curvature tensor $R$ with the curvature-like tensor $G$, $R = c\, G$}, similarly, \textit{a Riemannian manifold $(M^{n},g)$ is algebraically characterised to be a Deszcz symmetric space $M^{n}((L))$ by the functional proportionality of its curvature tensor $R\cdot R$ with its Tachibana tensor $\wedge_{g}\cdot R$, $R\cdot R = L\, \wedge_{g}\cdot R$ for some function $L:\, M^{n}\rightarrow \mathbb{R}$}.

As far as we know, it was via \textit{the curvature condition} $R\cdot R = L\, \wedge_{g}\cdot R$, $L:\, M^{n}\rightarrow\mathbb{R}$, that the pseudo symmetry of (semi) Riemannian spaces initially started to be studied. This condition as such appeared e.g. in studies by Wieslaw Grycak on semi symmetric warped product spaces and e.g. in studies by Sinyukov and Mike\v{s} and Venzi concerning geodesic mappings on semi symmetric spaces, while \textit{the term `'pseudo symmetric spaces''} for the (semi) Riemannian manifolds which satisfy this curvature condition, as far as we know, did first appear in an article by Deszcz and Grycak. In any case, in the 1970ties and 1980ties, \textit{the relevance of this intrinsically pseudo symmetric spaces became more clear through some investigations in the geometry of submanifolds, beginning with the studies by Ryszard Deszcz on the extrinsic spheres in the Szab\'{o} symmetric spaces}, as extensions of the studies of B.-Y. Chen and by Z. Olszak on the totally umbilical submanifolds of the Cartan symmetric spaces. On the origins and first publications concerning what later became to be known as \textit{the Deszcz symmetric spaces}, and also on related curvature conditions involving various \textit{other curvature tensors}, and also on the r\^{o}les played by pseudo symmetry in \textit{the theory of general relativity}, and for a first announcement about the study of \textit{extrinsically pseudo symmetric or pseudo parallel submanifolds}, which later a.o. resulted in the papers \cite{dillen,verstraelen1}, we refer to \cite{deszcz1,verstraelen2,haesen3,deszcz2}.


\section{Definition of quasi Einstein spaces}

\textit{The $(0,2)$ Ricci curvature tensor} $S$ of a Riemannian manifold $(M^{n},g)$, $n\geq 3$, is defined by $S(\tilde{X},\tilde{Y})=\sum_{t}R(\tilde{X},\tilde{E}_{t},\tilde{E}_{t},\tilde{Y})$, whereby $\tilde{X}$ and $\tilde{Y}$ are arbitrary tangent vector fields and $(\tilde{\varepsilon})=\left\{\tilde{E}_{1},\tilde{E}_{2},\ldots,\tilde{E}_{n}\right\}$, $(t,s\in\left\{1,2,\ldots,n\right\})$, is any orthonormal tangent frame field on $(M^{n},g)$, and \textit{the $(1,1)$ Ricci tensor} $S$ is defined by $g(S(\tilde{X}),\tilde{Y}) = S(\tilde{X},\tilde{Y})$. \textit{For a unit tangent vector $\tilde{U}$, the scalar} 

\[ \rho(\tilde{U}) = S(\tilde{U},\tilde{U}) = \sum_{t}R(\tilde{U},\tilde{E}_{t},\tilde{E}_{t},\tilde{U}) = \sum_{t} K(\tilde{U}\wedge\tilde{E}_{t}) \]

{\noindent}\textit{is called the Ricci curvature of $(M^{n},g)$ in the direction $\tilde{U}$}. The critical values of $\rho(\tilde{U})$ and the directions in which they are attained are called \textit{the Ricci principal curvatures or the principal Ricci curvatures} and \textit{the Ricci principal directions or the principal Ricci directions of $(M^{n},g)$}; it are the eigenvalues and the eigendirections of the Ricci tensor $S$.

\textit{Einstein spaces are the Riemannian manifolds $(M^{n},g)$, $n\geq 3$, of which the $(0,2)$ Ricci tensor $S$ is proportional to the metric tensor $g$}, $S=\rho\, g$, or, still, for which \textit{at every point the Ricci curvature function is isotropic}, that is, for which, \textit{the Ricci curvature $\rho(\tilde{U})$ is the same in all directions $\tilde{U}$}, or, still, \textit{Einstein spaces are the Riemannian manifolds which at every point have only one principal Ricci curvature $\rho$}, this $\rho$ then having \textit{multiplicity $n$} as sole eigenvalue of $S$ at the concerned point. As is well known: \textit{for all Einstein spaces $(M^{n},g)$, $n\geq 3$, this unique Ricci principal curvature is constant on $M^{n}$}, (cfr. \cite{kuhnel}). And, sometimes, for short, Einstein spaces may be called E spaces and be denoted as such.

\textit{(q) Quasi Einstein spaces}, or, for short, \textit{(q) QE spaces}, are defined as the Riemannian manifolds $(M^{n},g)$, $n\geq 3$, with \textit{precisely two distinct principal Ricci curvatures}, say $\rho$ and $\rho^{\perp}\neq\rho$, \textit{with fixed multiplicities}, say $q\geq 1$ and $q^{\perp}=n-q\geq q$. Such spaces do admit \textit{two differentiable and orthogonally complementary distributions}, say $D$ and $D^{\perp}$, \textit{with fixed dimensions} $q$ and $q^{\perp}$, namely, \textit{the eigenspaces of the Ricci tensor corresponding to its two distinct eigenvalues} $\rho$ and $\rho^{\perp}$. And, furtheron, (q) quasi Einstein spaces or (q) QE spaces may sometimes simply be called \textit{quasi Einstein spaces} or \textit{quasi E spaces}, and, then the actual multiplicities $q$ and $q^{\perp}$ of the two distinct Ricci principal curvatures $\rho$ and $\rho^{\perp}$, or, equivalently, the actual dimensions $q$ and $q^{\perp}$ of the two corresponding Ricci distributions $D$ and $D^{\perp}$, will be specified only when needed.

As far as we know, in case $q=1$, (q) quasi Einstein spaces have been studied for a long time and by many people and in these studies (1) quasi Einstein spaces were called `'quasi or pseudo or so Einstein spaces'' and in case $q>1$, as far as we know (q) quasi Einstein spaces so far did only occur in the literature sporadically and then they were called `'half and half Einstein spaces'' or `'partially Einstein spaces''.

At this stage, the following small aside might not be too much out of place here, we hope. Before the work of Ricci on his tensor and on his principal curvatures and principal directions, a significant first step on this way was made by Souvorov in his 1871 master thesis at Kazan, shortly after the first publications on Riemannian geometry, by Riemann and Helmholtz. \textit{Souvorov showed that at any point of a 3-dimensional Riemannian space $(M^{3},g)$ the critical values of the Riemannian sectional curvatures there are attained in three mutually orthogonal 2-planes}; (in retrospect, one can hardly avoid herein to see the origin of the theory on the Ricci curvatures for arbitrary Riemannian manifolds). And, as Marcel Berger wrote in his contribution to the book `'Chern - A Great Geometer of the Twentieth Century'': `'\textit{... knowing $K$ is knowing $R$ (...). But the relations between $K$ and $R$ are subtle and still not completely understood (e.g.  what the critical planes of $K$ are, and how they are distributed})''. Having originated in his answer to Shiing-Shen Chern's 1968 Kansas Lecture Notes' question \textit{to determine intrinsic geometric conditions on Riemannian manifolds $(M^{n},g)$ that would prevent their minimal isometric immersibility in Euclidean spaces $\mathbb{E}^{n+m}$, (with arbitrary codimension $m$}; and the only known such condition at that time being to have a non-negative definite Ricci tensor), \textit{Bang-Yen Chen's $\delta$-curvatures theory \cite{chen3} moreover has indeed been effectively contributing so much to this understanding the lack of which had been drawn attention to by Berger}. And, concerning the hereby occuring interplay between the extrinsic and the intrinsic geometries of submanifolds, here we confine to recall that the $\delta(2)$ Chen ideal submanifolds $M^{n}$ in $\mathbb{E}^{n+m}$ precisely do assume the very particular shapes for which the corresponding surface tension is as small as possible, that their mean curvature vector field does determine a first principal Casorati normal vector field and that their intrinsic principal Ricci directions do co-incide with their extrinsic tangent principal Jordan directions \cite{decu1,decu2}.

From Thurston's \cite{thurston} comes the following quote: \textit{`'What is geometry? Up till now, we have discussed three kinds of three-dimensional geometry: hyperbolic, Euclidean and spherical. They have in common the property of being as uniform as possible: their isometries can move any point to any other point (homogeneity), and can take any orthonormal frame in the tangent space at a point to any other orthonormal frame at that point (isotropy). There are more possibilities if we remove the isotropy condition, allowing the spaces to have a grain, so to speak, so that certain directions are geometrically distinguished from others}.

\textit{An enumeration of additional three-dimensional geometries depends on what spaces we wish to consider and what structures we use to define and distinguish the spaces. For instance, do we think of a geometry as a space equiped with such notions as lines and planes, or as a space equiped with a notion of congruence, or as a space equiped with either a metric or a Riemannian metric? (...) For logical purposes, we must pick only one definition. We choose to represent a geometry as a space equiped with a group of congruences, that is, a $(G,X)$-space.}

\textbf{Definition} A model geometry $(G,X)$ is a manifold $X$ together with a Lie group $G$ of diffeomorphisms of $X$, such that: (a) $X$ is connected and simply connected; (b) $G$ acts transitively on $X$, with compact point stabilizers; (c) $G$ is not contained in any larger group of diffeomorphisms of $X$ with compact stabilizers of points; (d) there exists at least one compact manifold modeled on $(G,X)$''.

And, \textit{in the following Theorem 3.8.4, then Thurston lists up his eight 3-dimensional model geometries}: the \textit{Euclidean geometry} $\mathbb{E}^{3}$ and the classical non-Euclidean \textit{spherical and hyperbolical geometries} $\mathbb{S}^{3}$ and $\mathbb{H}^{3}$, (which, in this order, are \textit{the model real space forms $M^{3}(0)$, $M^{3}(1)$ and $M^{3}(-1)$}, respectively), and, then \textit{the five Thurston model geometries which are non-isotropic, namely: $\mathbb{S}^{2}\times\mathbb{E}^{1}$, $\mathbb{H}^{2}\times\mathbb{E}^{1}$, $\tilde{SL}(2,\mathbb{R})$, $H_{3}$ (nilgeometry) and $Sol$ (solvgeometry)}.

\textit{Thurston's non-isotropic 3-dimensional model geometries are quasi Einstein spaces}: they have two distinct principal Ricci curvatures $\rho$ and $\rho^{\perp}$, with respective multiplicities 1 and 2; so, at everyone of their points, their tangent spaces essentially split up in the 1-dimensional and in the 2-dimensional mutually orthogonal eigenspaces $D$ and $D^{\perp}$ of their Ricci tensor, and, for tangent directions moving so to say from $D$ to $D^{\perp}$ their Ricci curvatures accordingly change nicely from $\rho$ to $\rho^{\perp}$. In some way, one could consider this situation as a kind of `'\textit{mild anisotropy}''. On the other hand, for the Riemannian spaces $(M^{3},g)$ of which the Ricci tensor has three mutually distinct principal curvatures $\rho_{1}$ and $\rho_{2}$ and $\rho_{3}$, each having multiplicity 1, at every point of $(M^{3},g)$, depending on their position in the 3-dimensional tangent space there, all different tangent directions basically do have different Ricci curvatures. From this point of view, such generic Riemannian spaces $(M^{3},g)$ are `'\textit{wildly anisotropic}'', in that, roughly speaking, all their tangent directions are geometrically distinguished from all the other tangent directions. And, so, also from this point of view, \textit{Thurston's extension from the real space forms $\mathbb{E}^{3}$ and $\mathbb{S}^{3}$ and $\mathbb{H}^{3}$ to his eight 3-dimensional model geometries by removing  the real space forms' isotropy quality does very well realise this goal}.

In Thurston's words `'\textit{when thinking of a geometry as a space equiped with a Riemannian metric'', the main isometrically invariant geometrically defined scalar value associated with a tangent direction of a Riemannian manifold $(M^{n},g)$ is its Ricci curvature}. And, \textit{the Riemannian spaces which are isotropic in the sense of Riemannian geometry}, that is, the Riemannian spaces for which all tangent directions do have the same Ricci curvature, \textit{are the Einstein spaces}. Next, taking into consideration that \textit{the real space forms in a way are  `' the most perfect Einstein spaces''}, one may well imagine that, for some geometers when wanting \textit{to remove from the real space forms their isotropy condition in the above sense of Riemannian geometry}, in their mind, \textit{a most natural option might be to replace it by the condition for Riemannian manifolds $(M^{n},g)$, $n\geq 3$, to be quasi Einstein spaces}.

The folllowing is a \textit{classical result of Schouten and Struik}, (cfr. \cite{kuhnel}).

\begin{mythm}{E} \label{thE}
A 3-dimensional Riemannian manifold is an Einstein space if and only if it is a real space form.
\end{mythm}

{\noindent}When \textit{restricting to the Riemannian manifolds $(M^{3},g)$ on which the multiplicities of the principal Ricci curvatures are constant}, in \cite{deszcz3,deszcz4} this result was extended as follows.

\begin{mythm}{F} \label{thF}
A 3-dimensional Riemannian manifold is Deszcz symmetric if and only if it is an Einstein space or it is a quasi Einstein space.
\end{mythm}

{\noindent}Thus, in particular, \textit{all eight 3-dimensional Thurston model geometric spaces} are Deszcz symmetric spaces, and, as a matter of fact, it all \textit{are Deszcz symmetric spaces of constant type}: more specifically, \textit{$\mathbb{E}^{3}$ and $\mathbb{S}^{3}$ and $\mathbb{H}^{3}$ and $\mathbb{S}^{2}\times\mathbb{E}^{1}$ and $\mathbb{H}^{2}\times\mathbb{E}^{1}$} are Cartan symmetric spaces, and so, as particular Szab\'{o} symmetric spaces, it \textit{are Deszcz symmetric spaces $M^{3}((0))$}; on the other hand, $\tilde{SL}(2,\mathbb{R})$ and $H_{3}$ are Deszcz symmetric spaces $M^{3}((I))$ and $Sol$ is a Deszcz symmetric space $M^{3}((-I))$, \cite{belkhelfa,trans2}, (whereby $I$ and $-I$ denote the functions $I:\, M^{n}\rightarrow\mathbb{R}:\, p\mapsto I(p)=1$ and $-I:\, M^{n}\rightarrow\mathbb{R}:\, (-I)(p)=-1$, respectively).

However, when going through \textit{the list of the nineteen or so 4-dimensional Thurston model geometries}, (for instance in \cite{wall}), then, \textit{beyond the Riemannian isotropic real space forms $\mathbb{E}^{4}$ and $\mathbb{S}^{4}$ and $\mathbb{H}^{4}$}, further one may encounter \textit{model spaces which are Riemannian anisotropic of all kinds}, going from the mildest possible anisotropy - quasi Einstein spaces - to the wildest possible anisotropy - Riemannian spaces of dimension 4 with 4 mutually distinct Ricci principal curvatures -.

Next follows another small aside, this one about two of the \textit{multiple uses of the term `'isotropy'' in geometry. In Riemannian geometry}, isotropy of 1D directions pretty naturally refers to the property that \textit{at all points the Ricci curvature is the same in all such directions}. \textit{In the above quotation from Thurston} concerning his model geometries, the `'\textit{isotropy of orthonormal tangent frames}'' together with the there mentioned `'\textit{1-point-homogeneity}'' may readily be thought of in connection with `'\textit{the axiom of free mobility}'', going back to Helmholtz and Riemann. Almost a century later, Jacques Tits may well have arrived at the latter's `'ultimately abstract version'' when \textit{characterising the real space forms as the `'3-point homogeneous spaces}'', while, still going the step further to \textit{the `'2-point homogeneous spaces}'', one basically arrived at \textit{the rank 1 symmetric spaces of Hsien-Chung Wang}, (cfr. the section `'The Space Problems'' in \cite{freudenthal}, and, in addition to the references given therein, see also \cite{szabo3}). 


\section{On the conformally Euclidean Deszcz symmetric spaces}

\textit{The studies of Gauss on the extension of the cartographers' stereographic and Mercator projections from round spheres on Euclidean planes to infinitesimally conformal maps between any two surfaces in $\mathbb{E}^{3}$}, were at the origin of \textit{complex analysis} and brought along \textit{the proof of the theorema egregium in isothermal co-ordinates}. And, in the context of these studies, Gauss did express his opnion that, \textit{for spaces to resemble each other well, the essential condition is that these spaces be similar in their smallest parts}, (cfr. \cite{dombrowski}). This condition may equivalently be formulated as \textit{that these spaces be locally conformal to each other}, or, for short, as \textit{that these spaces be conformal to each other}. Besides \textit{the natural ability in human vision to readily recognise similar things} despite their eventually rather different actual sizes and the fact that \textit{Euclidean geometry did originate as the science of human vision}, (cfr.\cite{klein,verstraelen3}), the above may help to see \textit{the significance of the class of Riemannian manifolds which are conformal to Euclidean spaces}, or, still, the class of the \textit{conformally Euclidean spaces}.

In this respect, for Riemannian manifolds $(M^{n},g)$, $n\geq 3$, we recall the definition of \textit{Weyl's $(0,4)$ conformal curvature tensor $C$}:

\newcommand{\tv}{\tilde{V}}

\begin{eqnarray}
C(\tilde{V}_{1},\tv_{2},\tv_{3},\tv_{4}) & = & R(\tv_{1},\tv_{2},\tv_{3},\tv_{4}) - \{ g(\tv_{1},\tv_{4})S(\tv_{2},\tv_{3}) - g(\tv_{1},\tv_{3})S(\tv_{2},\tv_{4}) \nonumber \\
 & & +g(\tv_{2},\tv_{3})S(\tv_{1},\tv_{4}) - g(\tv_{2},\tv_{4})S(\tv_{1},\tv_{3}) \}/(n-2)  \label{eqC} \\
 & & +\tau \{ g(\tv_{1},\tv_{4})g(\tv_{2},\tv_{3}) - g(\tv_{1},\tv_{3})g(\tv_{2},\tv_{4}) \}/(n-1)(n-2)\ , \nonumber 
\end{eqnarray}

{\noindent}whereby 

\[ \tau = \sum_{t}S(\tilde{E}_{t},\tilde{E}_{t}) = \sum_{t,s}R(\tilde{E}_{t},\tilde{E}_{s},\tilde{E}_{s},\tilde{E}_{t}) = \sum_{t,s}K(\tilde{E}_{t}\wedge \tilde{E}_{s})\ ,\ \ t\neq s\ , \]

{\noindent}is \textit{the scalar curvature} of $(M^{n},g)$. And the following \textit{classical results} are \textit{due to Jan Schouten and Hermann Weyl}, (cfr. \cite{kuhnel}).

\begin{mythm}{G} \label{thG}
For every 3-dimensional Riemannian manifold $C=0$, just like that.
\end{mythm}

\begin{mythm}{H} \label{thH}
For a Riemannian manifold $(M^{n},g)$ of dimension $n\geq 4$, $C=0$ if and only if $(M^{n},g)$ is a conformally Euclidean space.
\end{mythm}

In relation with Theorems \ref{thG} and \ref{thH}, and when \textit{restricting to Riemannian manifolds on which the multiplicities of the principal Ricci curvatures are constant}, like in Theorem \ref{thF}, from \cite{deprez} we recall the following.

\begin{mythm}{I} \label{thI}
A conformally Euclidean Riemannian manifold $(M^{n},g)$, $n\geq 4$, is Deszcz symmetric if and only if it is an Einstein space or it is a quasi Einstein space.
\end{mythm}

{\noindent}And in relation with Theorem \ref{thE}, the following \textit{classical result} is \textit{due to Jan Schouten and Dirk Jan Struik}, (cfr. \cite{kuhnel}).

\begin{mythm}{J}
A Riemannian manifold $(M^{n},g)$, $n\geq 4$, is a conformally Euclidean Einstein space if and only if it is a real space form.
\end{mythm}

\vspace{4mm}

The \textit{new results} of the present paper are the following.

\begin{theorem} \label{th1}
A Riemannian manifold of dimension $n\geq 4$ is a conformally Euclidean quasi Einstein space if and only if it is a quasi space form.
\end{theorem}

\begin{theorem} \label{th2}
A 3-dimensional Riemannian manifold is a quasi Einstein space if and only if it is a quasi space form.
\end{theorem}

\begin{theorem} \label{th3}
A 3-dimensional Riemannian manifold is Deszcz symmetric if and only if it is a real space form or it is a quasi space form.
\end{theorem}

\begin{theorem} \label{th4}
A conformally Euclidean Riemannian manifold of dimension $n\geq 4$ is Deszcz symmetric if and only if it is a real space form or it is a quasi space form.
\end{theorem}

\vspace{4mm}

{\noindent}\underline{\textit{Proofs}}

\vspace{2mm}

{\noindent}Clearly it suffices to prove the following statement: \textit{for a Riemannian manifold $(M^{n},g)$ of dimension $n\geq 3$, $C=0$ and $(M^{n},g)$ is a (q) quasi Einstein space if and only if $(M^{n},g)$ is a (q) quasi constant curvature space}.

First of all we take note of the fact that, as readily follows from (\ref{eqC}), \textit{on all Riemannian manifolds $(M^{n},g)$, $n\geq 3$, the sectional curvatures $K_{C}(\tilde{\pi})$ and $K(\tilde{\pi})$ of the curvature tensors $C$ and $R$ for tangent 2-planes spanned by orthonormal vectors $\tilde{X}$ and $\tilde{Y}$, $\tilde{\pi}=\tilde{X}\wedge\tilde{Y}$, are related as follows}, (cfr. also \cite{chen4}):

\begin{equation}
K_{C}(\tilde{X}\wedge\tilde{Y}) = K(\tilde{X}\wedge\tilde{Y}) - \left\{ \rho(\tilde{X}) + \rho(\tilde{Y})\right\}/(n-2) + \tau/(n-1)(n-2)\ . \label{eqKC}
\end{equation}

{\noindent}\textit{Now, we assume that $(M^{n},g)$, $n\geq 3$, is a (q) QE space for which $C=0$}. In particular, \textit{$(M^{n},g)$ being a (q) QE space}, let $(\varepsilon,\varepsilon^{\perp}) = \{E_{1},\ldots,E_{q},E_{1^{\perp}}^{\perp},\ldots,E_{(q^{\perp})^{\perp}}^{\perp}\}$ be any \textit{orthonormal tangent frame consisting of principal Ricci vectors} such that the $q$-dimensional and $q^{\perp}$-dimensional \textit{Ricci eigenspaces} $D$ and $D^{\perp}$ corresponding with the two \textit{distinct principal Ricci curvatures} $\rho$ and $\rho^{\perp}$ are given by $D=E_{1}\wedge\ldots\wedge E_{q}$ and $D^{\perp}=E_{1^{\perp}}^{\perp}\wedge\ldots\wedge E_{(q^{\perp})^{\perp}}^{\perp}$, respectively, or, still, such that $S(E_{i},E_{i})=\rho\neq \rho^{\perp} = S(E_{i^{\perp}}^{\perp},E_{i^{\perp}}^{\perp})$, $S(E_{i}) = \rho\, E_{i}$, $S(E_{i^{\perp}}^{\perp}) = \rho^{\perp}\, E_{i^{\perp}}^{\perp}$, $S(E_{i},E_{j}) = S(E_{i^{\perp}}^{\perp},E_{j^{\perp}}^{\perp}) = S(E_{i},E_{i^{\perp}}^{\perp})=0$, whereby $i,j\in\{1,2,\ldots,q\}$ and $i^{\perp},j^{\perp}\in\{1^{\perp},\ldots,(q^{\perp})^{\perp}\}$ and $i\neq j$, $i^{\perp}\neq j^{\perp}$. To the following we will hereafter refer to as to the Lemma.

\begin{lemma}
When applying (\ref{eqC}) to any four vectors from $(\varepsilon,\varepsilon^{\perp})$ of which at least three of these four are mutually distinct, then from $C=0$ it follows that $R(\, .\, ,\, .\, ,\, .\, ,\, .\, ) = 0$ for such four vectors.
\end{lemma}

Next, \textit{in case $q=1$}, by (\ref{eqKC}) it follows that

\begin{equation}
\forall i,i^{\perp}\ :\\ K(E_{i}\wedge E_{i^{\perp}}^{\perp}) = \{\rho+\rho^{\perp}\}/(n-2) - \tau\ . \label{eq3} 
\end{equation}

{\noindent}And thus, for any unit vector $X\in D$ and any unit vector $X^{\perp}\in D^{\perp}$, denoting their tangent 2-plane by $\overline{\pi}$, $\overline{\pi} = X\wedge X^{\perp}$, and, denoting the Riemannian sectional curvature of this plane by $\overline{K}$,

\begin{equation}
\overline{K} = K(\overline{\pi}) = \{\rho+\rho^{\perp}\}/(n-2) - \tau\ . \label{eq4}
\end{equation}

{\noindent}Similarly, by (\ref{eqKC}) it follows that

\begin{equation}
\forall i^{\perp}\neq j^{\perp}\ :\ \ K(E_{i^{\perp}}^{\perp}\wedge E_{j^{\perp}}^{\perp}) = \{\rho^{\perp}+\rho^{\perp}\}/(n-2) - \tau\ . \label{eq5}
\end{equation}

{\noindent}And thus, for any orthonormal vectors $X^{\perp},Y^{\perp}\in D^{\perp}$, denoting their tangent 2-plane by $\pi^{\perp}$, $\pi^{\perp}=X^{\perp}\wedge Y^{\perp}$, and, denoting the Riemannian sectional curvature of this plane $\pi^{\perp}$ by $K^{\perp}$, and, also since $\rho\neq \rho^{\perp}$,

\begin{equation}
K^{\perp} = K(\pi^{\perp}) = \{ 2 \rho^{\perp}\}/(n-2) - \tau\neq \overline{K}\ . \label{eq6}
\end{equation}

{\noindent}And so, for any tangent 2-plane $\tilde{\pi} = (X\cos\theta+X^{\perp}\sin\theta)\wedge Y^{\perp}$, whereby $X\in D$ and $X^{\perp},Y^{\perp}\in D^{\perp}$ are orthonormal and whereby $\theta\in \left[0,\frac{\pi}{2}\right]$, making use of the above Lemma and of (\ref{eq4}) and (\ref{eq6}), we get that

\begin{equation}
K(\tilde{\pi}) = \overline{K}\cos^{2}\theta + K^{\perp}\sin^{2}\theta\ , \label{eqKpi}
\end{equation}

{\noindent}\textit{which shows that $(M^{n},g)$ is a (1) QCC space}.

And next, \textit{in case $q>1$}, in addition to (\ref{eq4}) and (\ref{eq6}), by (\ref{eqKC}) it moreover follows that

\begin{equation}
\forall i\neq j\ :\ \ K(E_{i}\wedge E_{j}) = \{\rho+\rho\}/(n-2) - \tau\ . \label{eq7}
\end{equation}

{\noindent}And thus, for any orthonormal $X,Y\in D$, denoting their tangent 2-plane by $\pi$, $\pi=X\wedge Y$, and, denoting the Riemannian sectional curvature of this plane $\pi$ by $K$, and, also since $\rho\neq \rho^{\perp}$,

\begin{equation}
K = K(\pi) = \{2\rho\}/(n-2) - \tau\neq K^{\perp}\ . \label{eq8}
\end{equation}

{\noindent}And so, for any tangent 2-plane $\tilde{\pi} = (X\cos\theta+X^{\perp}\sin\theta)\wedge(Y\cos\varphi\wedge Y^{\perp}\sin\varphi)$, whereby $X,Y\in D$ and $X^{\perp},Y^{\perp}\in D^{\perp}$ are orthonormal and whereby $\theta,\varphi\in\left[0,\pi/2\right]$, making use of the above Lemma and of (\ref{eq4}) and (\ref{eq6}) and (\ref{eq8}), we get that

\begin{equation}
K(\tilde{\pi}) = K\cos^{2}\theta\cos^{2}\varphi + K^{\perp}\sin^{2}\theta\sin^{2}\varphi +\overline{K}(\cos^{2}\theta\sin^{2}\varphi+\sin^{2}\theta\cos^{2}\varphi)\ , \ \ \overline{K} = (K+K^{\perp})/2\ , \label{eqstarstar}
\end{equation}

{\noindent}\textit{which shows that $(M^{n},g)$ is a (q) QCC space}.

\textit{Conversely, now we assume that $(M^{n},g)$, $n\geq 3$, is a (q) QCC space}. To begin with, $D$ and $D^{\perp}$ being the orthogonally complementary $q$- and $q^{\perp}=(n-q)$-dimensional distributions of the given (q) QCC space, let $(\varepsilon,\varepsilon^{\perp})=\{ E_{1},\ldots,E_{q},E_{1^{\perp}}^{\perp},\ldots,E_{(q^{\perp})^{\perp}}^{\perp}\}$ be an orthonormal tangent frame such that $\forall\, i:\ E_{i}\in D$ and $\forall\, i^{\perp}:\ E_{i^{\perp}}^{\perp}\in D^{\perp}$. Then, making use of (\ref{eqKpi}) and (\ref{eqstarstar}), we calculate \textit{the Ricci curvatures of the vectors $E_{i}$ and $E_{i^{\perp}}^{\perp}$}:

\begin{eqnarray}
\rho(E_{i}) & = & \sum_{j\neq i} K(E_{i}\wedge E_{j}) + \sum_{i^{\perp}} K(E_{i}\wedge E_{i^{\perp}}^{\perp}) \nonumber \\
 & = & (q-1) K + q^{\perp} \overline{K}\ , \label{eq9}
\end{eqnarray}

\begin{eqnarray}
\rho(E_{i^{\perp}}^{\perp}) & = & \sum_{i} K(E_{i^{\perp}}^{\perp}\wedge E_{i}) + \sum_{j^{\perp}\neq i^{\perp}} K(E_{i^{\perp}}^{\perp}\wedge E_{j^{\perp}}^{\perp}) \nonumber \\
 & = & q \overline{K}+(q^{\perp}-1)K^{\perp}\ . \label{eq10}
\end{eqnarray}

{\noindent}Thus, for all unit vectors $U\in D$ and $U^{\perp}\in D^{\perp}$:

\begin{equation}
\rho(U) = \rho\neq \rho^{\perp} = \rho(U^{\perp})\ , \label{eq11}
\end{equation}

{\noindent}while the precise relations between these Ricci curvatures $\rho$ and $\rho^{\perp}$ in the distributions $D$ and $D^{\perp}$ of a (q) QCC space and its sectional curvatures $\overline{K}$ and $K^{\perp}$ and $K$ are given in (\ref{eq9}) and (\ref{eq10}). Next we calculate the Ricci curvature $\rho(\psi) = S(\tilde{U},\tilde{U})$ in the direction of any unit tangent vector $\tilde{U}= U\, \cos\psi+U^{\perp}\, \sin\psi$ whereby $\psi\in\left[0,\pi/2 \right]$. In order to do so, we consider an orthonormal frame $(\tilde{\varepsilon},\tilde{\varepsilon}^{\perp}) = \{ \tilde{U} = U\, \cos\psi+U^{\perp}\, \sin\psi=\tilde{U}_{1}, \tilde{U}^{\perp}=-U\sin\psi+U^{\perp}\cos\psi = \tilde{U}_{1^{\perp}}^{\perp},U_{2},\ldots,U_{q}, U_{2^{\perp}}^{\perp},\ldots, U_{(q^{\perp})^{\perp}}^{\perp}\}$, whereby $U_{2},\ldots,U_{q}\in D$ are perpendicular to $U$ and whereby $U_{2^{\perp}}^{\perp},\ldots, U_{(q^{\perp})^{\perp}}^{\perp}\in D^{\perp}$ are perpendicular to $U^{\perp}$. Then,

\[ \rho(\psi) = \rho(\tilde{U}) = K(\tilde{U}\wedge\tilde{U}^{\perp}) + \sum_{i\neq 1} K(\tilde{U}\wedge U_{i}) + \sum_{i^{\perp}\neq 1^{\perp}} K(\tilde{U}\wedge U_{i^{\perp}}^{\perp})\ . \]

{\noindent}In case $q=1$, using (\ref{eqKpi}) this gives that

\begin{equation}
\rho(\psi) = \overline{K}+(n-2)(\overline{K}\cos^{2}\psi + K^{\perp}\sin^{2}\psi)\ , \label{eq11}
\end{equation}

{\noindent}and, in case $q>1$, using (\ref{eqstarstar}) this gives that

\begin{equation}
\rho(\psi) = \overline{K} +(q-1)(K\cos^{2}\psi+\overline{K}\sin^{2}\psi) + (n-q-1)(K^{\perp}\sin^{2}\psi+\overline{K}\cos^{2}\psi)\ . \label{eq12}
\end{equation}

{\noindent}From (\ref{eq11}) and (\ref{eq12}), respectively we get the following:

\begin{equation}
\frac{\mbox{d}\rho(\psi)}{\mbox{d}\psi} = 2(n-2)\cos\psi\sin\psi (K^{\perp}-\overline{K}) \label{eq13}
\end{equation}

{\noindent}and

\begin{equation}
\frac{\mbox{d}\rho(\psi)}{\mbox{d}\psi} = (n-2)\cos\psi\sin\psi (K^{\perp}-K)\ . \label{eq14}
\end{equation}

{\noindent}Given that on a (q) QCC space $K\neq K^{\perp}\neq \overline{K}$, (\ref{eq13}) and (\ref{eq14}) imply that in both cases, $q=1$ and $q>1$, \textit{the Ricci principal curvatures are attained in the directions $\tilde{U}=U\, \cos\psi+U^{\perp}\, \sin\psi$ for which $\psi=0$ or $\psi=\frac{\pi}{2}$}. \textit{This shows that a (q) QCC space is a (q) QE space}, with $\rho$ and $\rho^{\perp}\neq \rho$ as its two distinct principal Ricci curvatures and with $D$ and $D^{\perp}$ as its corresponding $q$- and $q^{\perp}=(n-q)$-dimensional principal Ricci distributions.

It further remains to be shown that $C=0$. And this can be done by straightforward verification as follows. On a (q) QCC space, for all their tangent 2-planes the Riemannian sectional curvatures are given by (\ref{eqKpi}) in case $q=1$ and by (\ref{eqstarstar}) in case $q>1$. By (\ref{eq9}) and (\ref{eq10}), rewritten hereafter as (\ref{eq9a}) and (\ref{eq10a}),

\begin{eqnarray}
\rho & = & (q-1)K+q^{\perp}\overline{K}\ , \label{eq9a} \\
\rho^{\perp} & = & q\overline{K} + (q^{\perp}-1)K^{\perp}\ , \label{eq10a}
\end{eqnarray}

{\noindent}the correspondences between $\rho$, $\rho^{\perp}$ on the one hand and $\overline{K}$, $K^{\perp}$ and $K$ on the other hand are clear, and, moreover we have that

\begin{equation}
\tau = \mbox{tr}\, S  =q\rho+q^{\perp}\rho^{\perp}\ . \label{eq15}
\end{equation}

{\noindent}Inserting these data into (\ref{eqKC}) for all possible tangent 2-planes $\tilde{\pi}$ of a (q) QCC space, it results that \textit{all the Weyl sectional curvatures $K_{C}(\tilde{\pi})$ do vanish, $\forall\, \tilde{\pi}:\ K_{C}(\tilde{\pi}) = 0$. And, since the conformal curvature tensor $C$ of Weyl is a generalised curvature tensor, this is equivalent with the vanishing of this curvature tensor itself}, and, thus, $C=0$, (cfr. \cite{kuhnel}). \hfill (Q.E.D.)

\vspace{4mm}

\textit{The real space forms of dimension $\geq 3$, or, still, the CC spaces $M^{n}(c)$, $n\geq 3$}, in a trivial way, \textit{are Deszcz symmetric spaces}, ($R\cdot R = L\, \wedge_{g}\cdot R$, $L:\, M^{n}\rightarrow\mathbb{R}$), \textit{with vanishing conformal curvature tensor of Weyl}, ($C=0$). In a way, the real space forms are the most symmetric of all Deszcz symmetric spaces, and, in particular, it already are very special Cartan symmetric spaces, ($\nabla R=0$), and, thus, even more special Szab\'{o} symmetric spaces, ($R\cdot R=0$). In any case, for the real space forms the double sectional curvature vanishes, $L=0$, or, still, \textit{the real space forms are very special Deszcz symmetric spaces $M^{n}((0))$}, (whereby $0:\, M^{n}\rightarrow\mathbb{R}\, :\, p\mapsto 0(p)=0$ is the null function on $M^{n}$). As asserted in Theorems \ref{th3} and \ref{th4}, \textit{besides the real space forms or the CC spaces, the quasi space forms, or, still, the QCC spaces $M^{n}(\overline{K},K^{\perp})$, $n\geq 3$, are the other Deszcz symmetric spaces with vanishing conformal curvature tensor of Weyl}. Among all Deszcz symmetric spaces, also \textit{the quasi space forms do form a particularly special subclass of extremely symmetric spaces, of which we now determine the double sectional curvature $L$}. In order to do so, let us consider two tangent 2-planes $\tilde{\pi}_{1}$ and $\tilde{\pi}_{2}$ given by $\tilde{\pi}_{1}=X\wedge X^{\perp}$ and $\tilde{\pi}_{2}=\tilde{X}\wedge Y^{\perp}$, whereby $X\in D$ and $X^{\perp},Y^{\perp}\in D^{\perp}$ are orthonormal, and, $D$ and $D^{\perp}$ are the two Ricci principal distributions of the quasi Einstein spaces $M^{n}(\overline{K},K^{\perp})$, in short notation: $S(D)=\rho\, D$ and $S(D^{\perp})=\rho^{\perp}D^{\perp}$, and whereby $\tilde{X}=(X+X^{\perp})/\sqrt{2}$. Then, a.o. making use of the above Lemma mentioned in the proofs, we find that

\begin{equation}
(\wedge_{g}\cdot R)(\tilde{X},Y^{\perp},Y^{\perp},\tilde{X};X,X^{\perp}) = K^{\perp} -\overline{K}\ , \label{eq16}
\end{equation}

{\noindent}(showing in particular, since $K^{\perp}\neq \overline{K}$, that \textit{such 2-planes $\tilde{\pi}_{1}$ and $\tilde{\pi}_{2}$ are curvature-dependent}), and that

\begin{equation}
(R\cdot R)(\tilde{X},Y^{\perp},Y^{\perp},\tilde{X};X,X^{\perp}) = \overline{K}\, (K^{\perp}-\overline{K})\ . \label{eq17}
\end{equation}

{\noindent}From (\ref{eq16}) and (\ref{eq17}) it follows that $L=\overline{K}$, which fact, for the record, will be formulated in the following.

\begin{proposition}
The QCC spaces $M^{n}(\overline{K},K^{\perp})$ are Deszcz symmetric spaces $M^{n}((\overline{K}))$.
\end{proposition}

\textit{In case $\overline{K} = K(\overline{\pi}) = K(X\wedge X^{\perp}) = 0$}, the (q) QCC spaces under consideration are Szab\'{o} symmetric, for either $q=1$ or $q>1$, whereby in the latter situation $K^{\perp}=-K$, and so, together with a later comment, one may recover the following \textit{classification result of K. Sekigawa and H. Takagi \cite{sekigawa}}.

\begin{mythm}{K} \label{thK}
The Szab\'{o} symmetric spaces $M^{n}((0))$, $n\geq 3$, with $C=0$ locally are (i) real space forms $M^{n}(c)$, or, (ii) product spaces $M^{q}(c)\times M^{n-q}(-c)$, $c\neq 0$, that is, products of real space forms with non-zero opposite curvatures, or, (iii) product spaces of real space forms $M^{n-1}(c)$, $c\neq 0$, with curves $M^{1}$.
\end{mythm}

The just mentioned comment refers to the following result from \cite{deprez}.

\begin{mythm}{L} \label{thL}
Let $(M^{n},g)$, $n\geq 3$, be a Riemannian manifold with vanishing conformal curvature tensor of Weyl. Then $(M^{n},g)$ is Deszcz symmetric, $R\cdot R=L\, \wedge_{g}\cdot R$, if and only if $(M^{n},g)$ is Ricci Deszcz symmetric, $R\cdot S = L\, \wedge_{g}\cdot S$.
\end{mythm}

{\noindent}A Riemannian manifold $(M^{n},g)$, $n\geq 3$, is said to be \textit{Ricci Deszcz symmetric} if it satisfies the curvature condition $R\cdot S = L\, \wedge_{g}\cdot S$ for some function $L:\, M^{n}\rightarrow\mathbb{R}$. Theorem \ref{thL} states that, whenever $C=0$, this condition is equivalent to the curvature condition $R\cdot R = L\, \wedge_{g}\cdot R$ to be Deszcz symmetric. And, in particular, Theorem \ref{thL} states that, whenever $C=0$, Riemannian manifolds $(M^{n},g)$, $n\geq 3$, are semi symmetric, $R\cdot R=0$, if and only if they are Ricci semi symmetric, $R\cdot S=0$, and the former Theorem \ref{thK} of Sekigawa-Takagi originally was formulated in this latter way. \textit{Geometrically the conditon $R\cdot S=0$ means the following} \cite{jahanara}: after parallel transports fully around infinitesimal co-ordinate parallelograms cornered at points $p$, tangent vectors $\vec{v}$ at $p$ in general are transformed to tangent vectors $\vec{v}^{\star}\neq\vec{v}$ at $p$, and $R\cdot S=0$ if and only if for all vectors $\vec{v}$ at all points $p$ after all such parallel transports $\rho(\vec{v}^{\star})=\rho(\vec{v})$, that is, if and only if, \textit{although in general tangent directions may change under such parallel transports, on Ricci semi symmetric spaces, their Ricci curvatures do remain invariant}.

For \textit{quasi Einstein spaces}, in passing, in Section 3 was mentioned \textit{the gradual change of the values of the Ricci curvatures from the value $\rho$ for directions in the distribution $D$ to the value $\rho^{\perp}\neq \rho$ for directions in the distribution $D^{\perp}$}. More precisely this means the following. Let $\tilde{X}=X\cos\theta+X^{\perp}\sin\theta$, $\theta\in \left[0,\pi/2\right]$, be an arbitrary tangent unit vector on a QE space, given in its canonical decomposition according to the direct sum $TM = D\oplus D^{\perp}$. Then, the Ricci curvature $\rho(\tilde{X})$ in the direction $\tilde{X}$ is given by

\begin{eqnarray*}
\rho(\theta) = \rho(\tilde{X}) = S(\tilde{X},\tilde{X}) & = & S(X\cos\theta+X^{\perp}\sin\theta,\ X\cos\theta+X^{\perp}\sin\theta) \\
 & = & S(X,X)\cos^{2}\theta + S(X^{\perp},X^{\perp})\sin^{2}\theta \\
 & = & \rho\, \cos^{2}\theta + \rho^{\perp}\, \sin^{2}\theta\ .
\end{eqnarray*}

{\noindent}We finish by referring to Bang-Yen Chen’s book ‘Geometry of Submanifolds’ \cite{chen73} of half a century ago for the start of the journey which is described in the present paper and to the Foreword of Bang-Yen Chen’s recent book \cite{chen17} on the geometry of warped products concerning the original rôle played herein by the (1)QCC spaces and for its discussion of Thurston’s geometrical space form problem.



\begin{thebibliography}{10}

\bibitem{kuhnel} W. K\"{u}hnel, \textit{Differentialgeometrie. Kurven - Fl\"{a}chen - Mannigfaltigheiten}, Vieweg, Wiesbaden (2008); English translation: \textit{Differential Geometry. Curves - Surfaces - Manifolds}, AMS Student Mathematical Library 16 (2006).

\bibitem{risteki} (a) I.B. Risteki and K.Tren\v{c}evski, \textit{Principal values and principal subspaces of two subspaces of vector spaces with inner product}, Beitr. Algebra Geom. \textbf{42} (2001), 289-300. \newline

(b) H. Gunawan, O. Neswan and W. Setya-Budhi, \textit{A formule for angles between subspaces of inner product spaces}, Beitr. Algebra Geom. \textit{46} (2009), 311-320.

\bibitem{chen1} B.-Y. Chen and K. Yano, \textit{Hypersurfaces of a conformally flat space}, Tensor N.S. \textbf{26} (1972), 318-322.

\bibitem{chen2} B.-Y. Chen and K. Yano, \textit{Special conformally flat spaces and canal hypersurfaces}, T\^{o}hoku Math. J. \textbf{25} (1973), 177-184.

\bibitem{yano1} K. Yano, C.-S. Houh and B.-Y. Chen, \textit{Intrinsic characterization of certain conformally flat spaces}, Kodai Math. Sem. Rep. \textit{25} (1973), 357-361.

\bibitem{boju} B. Boju and M. Popescu, \textit{Espaces \`{a} courbure quasi-constante}, J. Diff. Geom. \textbf{13} (1978), 373-383.

\bibitem{ganchev1} G. Ganchev and V. Mihova, \textit{Riemannian manifolds of quasi-constant sectional curvature}, J. reine angew. Math. \textbf{522} (2000), 119-141.

\bibitem{ganchev2} G. Ganchev and V. Mihova, \textit{A classification of Riemannian manifolds of quasi-constant sectional curvatures}, Ann. Sofia Univ. Fac. Math. and Inf. \textbf{102} (2015), 115-138.

\bibitem{trans} H. Weyl, \textit{Symmetry}, Princeton University Press, Princeton (1952).

\bibitem{trans2} S. Haesen and L. Verstraelen, \textit{Natural Intrinsic Geometrical Symmetries}, SIGMA \textbf{5} (2009), pp. 15, Special Issue `'\textit{\'{E}lie Cartan and Differential Geometry}''.

\bibitem{schouten} J.A. Schouten, \textit{Die direkte Analysis zur neueren Relativit\"{a}tstheorie}, Verhan. Konink. Akad. Wet. Amsterdam \textbf{12} (1918), nr. 6, 1-95.

\bibitem{matveev} V.S. Matveev, \textit{Geometric explanation of the Beltrami theorem}, Int. J. Methods Mod. Phys. \textbf{3} (2006), 623-629.

\bibitem{berger1} M. Berger, \textit{A Panoramic View of Riemannian Geometry}, Springer, Berlin (2003).
 
\bibitem{berger2} M. Berger, \textit{La g\'{e}om\'{e}trie m\'{e}trique de vari\'{e}t\'{e}s riemanniennes (...)}, in ``\textit{\'{E}lie Cartan et les math\'{e}matiques d' aujourd'hui}'', Ast\'{e}risque, Paris (1985), 9-66.

\bibitem{haesen2} S. Haesen and L. Verstraelen, \textit{Properties of a scalar curvature invariant depending on two planes}, Manuscripta Math. \textbf{122} (2007), 59-72.

\bibitem{szabo1} Z.I. Szab\'{o}, \textit{Structure theorems on Riemannian spaces satisfying $R(X,Y)\cdot R = 0$. I. The local version}, J. Diff. Geom. \textbf{17} (1982), 531-582.

\bibitem{szabo2} Z.I. Szab\'{o}, \textit{Structure theorems on Riemannian spaces satisfying $R(X,Y)\cdot R=0$. II. The global version}, Geom. Dedicata \textbf{19} (1985), 65-108.

\bibitem{levy} H. Levy, \textit{Tensors determined by a hypersurface in Riemannian space}, Trans. AMS \textbf{28} (1926), 671-694.

\bibitem{mikes1} J. Mike\v{s}, V. Kiosak and A. Van\v{z}urov\'{a}, \textit{Geodesic mappings of manifolds with affine connection}, Palack\'{y} University, Olomouc (2008).

\bibitem{mikes2} J. Mike\v{s}, A. Van\v{z}urov\'{a} and I. Hinterleitner, \textit{Geodesic mappings and some generalizations}, Palack\'{y} University, Olomouc (2009).

\bibitem{boeckx} E. Boeckx, O. Kowalski and L. Vanhecke, \textit{Riemannian manifolds of conullity two}, World Scientific, Singapore (1996).

\bibitem{lumiste} \"{U}. Lumiste, \textit{Semi-parallel submanifolds in real space forms}, Springer, Berlin (2009).

\bibitem{kowalski1} O. Kowalski and M. Sekizawa, \textit{Pseudo-symmetric spaces of constant type in dimension three - elliptic case}, Rend. Mat. Appl. (7) \textbf{17} (1997), 477-512.

\bibitem{kowalski2} O. Kowalski and M. Sekizawa, \textit{Pseudo-symmetric spaces of constant type in dimension three - non-elliptic case}, Bull. Tokyo Gakugei Univ. (4) \textbf{50} (1998), 1-28.

\bibitem{defever} F. Defever and R. Deszcz, \textit{A note on geodesic mappings of pseudosymmetric Riemannian manifolds}, Colloq. Math. \textbf{62} (1991), 313-319.

\bibitem{albujer} A. Albujer and S. Haesen, \textit{A geometrical interpretation of the null sectional curvature}, J. Geom. Phys. \textbf{60} (2010), 471-476.

\bibitem{eisenhart} L.P. Eisenhart, \textit{Riemannian Geometry}, Princeton University Press, Princeton (edition 1997).

\bibitem{dillen} F. Dillen, J. Fastenakels, S. Haesen, J. Van der Veken and L. Verstraelen, \textit{Submanifold theory and Parallel transport}, Kragujevac J. Math. \textbf{37} (2013), 33-43.

\bibitem{verstraelen1} L. Verstraelen, \textit{Natural extrinsic geometrical symmetries - an introduction}, AMS Contempory Mathematics, Vol. 674 (2016), 5-16; (from the `'\textit{AMS Special Session on Recent Advances in the Geometry of Submanifolds Dedicated to the Memory of Franki Dillen (1963-2013)}'', (eds. B. Suceava e.a.).

\bibitem{deszcz1} R. Deszcz, \textit{On pseudosymmetric spaces}, Bull. Soc. Math. Belg. A \textbf{44} (1992), 1-34.

\bibitem{verstraelen2} L. Verstraelen, \textit{Comments on the pseudo-symmetry in the sense of Ryszard Deszcz}, Geometry and Topology of Submanifolds VI, (eds. F. Dillen e.a.), World Scientific, Singapore (1994), 199-209.

\bibitem{haesen3} S. Haesen and L. Verstraelen, \textit{Chapter 8: Curvatures and Symmetries of Parallel Transport, and Chapter 9: Extrinsic Symmetries of Parallel transport}, in \textit{Differential Geometry and Topology, Discrete and Computational Geometry}, (eds. J.-M. Morvan and M. Boucetta), NATO Science Series - IOS Press (2005), 197-255.

\bibitem{deszcz2} R. Deszcz, S. Haesen and L. Verstraelen, \textit{Chapter 6: On Natural Symmetries}, in \textit{Topics in Differential Geometry}, (eds. A. \& I. Mihai and R. Miron), Ed. Acad. Rom\^{a}ne, Bucharest (2008), 249-308.

\bibitem{chen3} B.-Y. Chen, \textit{Pseudo-Riemannian Geometry, $\delta$-Invariants and Applications}, World Scientific, Singapore (2011).

\bibitem{decu1} S. Decu, A. Panti\'{c}, M. Petrovi\'{c}-Torg\v{a}sev and L. Verstraelen, \textit{Ricci and Casorati principal directions of $\delta(2)$ Chen ideal submanifolds}, Kragujevac Math. J. \textbf{37} (2013), 25-31.

\bibitem{decu2} S. Decu, B. Janahara, M. Petrovi\'{c}-Torg\v{a}sev and L. Verstraelen, \textit{On the Chen character of $\delta(2)$ ideal submanifolds}, Kragujevac Math. J. \textbf{32} (2009), 37-46.

\bibitem{thurston} W. Thurston, \textit{The Geometry and Topology of 3-Manifolds}, Lecture Notes, Princeton University Press, Princeton (1979).

\bibitem{deszcz3} R. Deszcz, L. Verstraelen and \c{S}. Yaprak, \textit{Pseudo-symmetric hypersurfaces in 4-dimensional spaces of constant curvature}, Bull. Inst. Math. Acad. Sinica \textbf{22} (1994), 167-179.

\bibitem{deszcz4} R. Deszcz, L. Verstraelen and \c{S}. Yaprak, \textit{Warped products realizing a certain condition of pseudosymmetry type imposed on the Weyl curvature tensor}, Chinese J. Math. \textbf{22} (1994), 139-157.

\bibitem{belkhelfa} M. Belkhelfa, R. Deszcz and L. Verstraelen, \textit{Symmetry properties of 3-dimensional d'Atri spaces}, Kyungpook Math. J. \textbf{46} (2006), 367-376.

\bibitem{wall} C. T. C. Wall, \textit{Geometries and Geometrical Structures in the real dimension 4 and the complex dimension 2}, Springer Lecture Notes 1167, eds. Dold and Eckmann; Proc. Univ. of Maryland Special Year of Low Dimensional Topology , College Park, 1983-1984.

\bibitem{freudenthal} H. Freudenthal and H.-G. Steiner, \textit{Chapter 13: Group Theory and Geometry}, in \textit{Fundamentals of Mathematics, Vol. II: Geometry}, (eds. H. Behnke e.a.), MIT Press, Cambridge (1974), 516-533.

\bibitem{szabo3} Z. I. Szab\'{o}, \textit{A short topological proof of the symmetry of 2-point homogeneous spaces}, Invent. math. \textbf{106} (1991), 61-64.

\bibitem{dombrowski} P. Dombrowski, \textit{150 years after Gauss' `'Disquisitiones generalis circa superficias curvas''}, Ast\'{e}risque, Paris (1979).

\bibitem{klein} F. Klein, \textit{Elementary mathematics from an advanced standpoint - Geometry}, Dover, New York (1939).

\bibitem{verstraelen3} L. Verstraelen, \textit{Submanifold theory - A contemplation of submanifolds}, AMS Contemporary Mathematics, Vol. 756 (2020), 21-56; (from the `'\textit{AMS Special Session on the Geometry of Submanifolds in Honor of Bang-Yen Chen}'', (eds. J. Van der Veken e.a.)).

\bibitem{deprez} J. Deprez, R. Deszcz and L. Verstraelen, \textit{Examples of pseudo-symmetric conformally flat warped products}, Chinese J. Math. \textbf{17} (1989), 51-65.

\bibitem{chen4} B.-Y. Chen, F. Dillen, L. Verstraelen and L. Vrancken, \textit{Characterizations of Riemannian space forms, Einstein spaces and conformally flat spaces}, Proc. AMS \textbf{128} (1999), 589-598.

\bibitem{sekigawa} K. Sekigawa and H. Takagi, \textit{On the conformally flat spaces satisfying a curvature condition on the Ricci tensor}, T\^{o}hoku Math. J. \textbf{23} (1971), 1-11.

\bibitem{jahanara} B. Jahanara, S. Haesen, Z. \c{S}ent\"{u}rk and L. Verstraelen, \textit{On the parallel transport of the Ricci curvatures}, J. Geom. Phys. \textbf{57} (2007), 1771-1777.

\bibitem{chen73} B.-Y. Chen, \textit{Geometry of Submanifolds}, Marcel Dekker, New York (1973).

\bibitem{chen17} B.-Y. Chen, \textit{Differential Geometry of Warped Product Manifolds and Submanifolds}, World Scientific, Singapore (2017).


\end{thebibliography}
\end{document}